\documentclass{amsart}
\usepackage[centertags]{amsmath}
\usepackage{amsfonts}
\usepackage{amssymb}
\usepackage{amsthm}
\usepackage{newlfont}

\newcommand{\R}{\mathbb R}
\newcommand{\N}{\mathbb N}
\newcommand{\1}{\mathbb I}
\newcommand{\Hs}{\mathfrak{H}}

\newcommand{\C}{\mathbb C}

\theoremstyle{plain}
\newtheorem{thm}{Theorem}[section]

\newtheorem{lem}[thm]{Lemma}
\newtheorem{prop}[thm]{Proposition}
\theoremstyle{definition}
\newtheorem{defn}{Definition}[section]
\theoremstyle{remark}\newtheorem{rem}{Remark}[section]
\newtheorem{ex}[thm]{Example}

\begin{document}

\title[Stochastic Banach Principle]
 {Stochastic Banach Principle in Operator Algebras}

\author{ Genady Ya. Grabarnik}

\address{Genady Ya. Grabarnik, IBM T.J. Watson Research Center, 19 Skyline dr, Hawthorne, 10532, USA }

\email{genady@us.ibm.com}

\author{Laura Shwartz}
\address{Laura Shwartz, Department of Mathematics, Applied Mathematics and
Astronomy, University of South Africa, Pretoria 0003, South Africa}
\email{lauralsh@gmail.com}



\date{January 24, 2006}

\dedicatory{}

\thanks{The second author is thankful to Dr. Louis E Labuschagne
(UNISA, South Africa) for constant support and carefull reading of the paper.}


\begin{abstract}
Classical Banach principle is an essential tool for the
investigation of the ergodic properties of $\check{C}$esaro
subsequences.

The aim of this work is to extend Banach principle to the case of
the stochastic convergence in the operator algebras.

We start by establishing a sufficient condition for the stochastic
convergence (stochastic Banach principle). Then we formulate
stochastic convergence for the bounded Besicovitch sequences, and,
as consequence for uniform subsequences.
\end{abstract}

\maketitle

\section{\protect\smallskip Introduction and Preliminaries}

In this paper we establish a Stochastic Banach Principle. The Banach
Principle is one of the most useful tools in "classical" point-wise
ergodic theory. The Banach principle was used to give an alternative proof
of the Birkhoff- Khinchin individual ergodic theorem. Typical
applications of the Banach Principle are Sato's theorem for
uniform subsequences \cite{Sato} and individual ergodic theorem for
the Besicovitch Bounded sequences \cite{Ryll} . Non-commutative
analogs for the (double side) almost everywhere convergence may be
found in papers \cite{GoldLit}, \cite{Chil}.

In this paper we establish a Banach Principle for convergence in
measure (Stochastic Banach Principle, Theorem 3.2.3).
We reformulate the theorem in a form convenient for applications (Theorem 3.2.4).
Based on the principle we
give a simplified proof of the stochastic ergodic theorem (compare
with \cite{GrabKatz1995}). We establish stochastic convergence for
Sato's uniform subsequences (Theorem 3.6) and a stochastic ergodic theorem for the
Besicovitch Bounded sequences (Theorem 3.7).

Note that these results are new even in the commutative case.

\par Throughout the  paper we denote by $M$ a von Neumann algebra with
semi-finite normal faithful trace $\tau$ acting on Hilbert space
$\mathfrak{H}$. Denote by $P(M)$ the set of all orthogonal projections
in $M$.

\par Recall the following definitions (combined from the papers by Segal
\cite{Segal}, Nelson \cite{Nel74}, Yeadon \cite{Yeadon76},  Fack and
Kosaki \cite{FackKos86} ):

\begin{defn}
A densely defined closed operator $x$ affiliated with von Neumanm
algebra $M$ is called \textbf{($\tau$) measurable} if for every
$\epsilon>0$ there exists projection $e\in P(M)$ with
$\tau(\1-e)<\epsilon$ such that $e(\mathfrak{H})\subset
\mathfrak{D}(x)$, where $\mathfrak{D}(x)$ is a domain of $x$.

Space of all \textbf{($\tau$) measurable operators affiliated with
$M$} is denoted by $S(M)$.

For convenience for a self-adjoint $x\in S(M)$ we denote by $\{x>t\}$
the \textbf{ spectral projection of $x$ corresponding to the interval
$(t,\infty]$}.
\end{defn}

\begin{defn}
\label{d:STOCH} Sequence $\{x_n\}_{n=1}^\infty$ \textbf{converges to
$0$ in measure (stochastically)} if for every $\epsilon>0$ and
$\delta>0$ there exists an integer $ N_0$  and a set of projections
$ \{e_n\}_{n\geq N_0} \subset P(M) $ such that $\| x_n\, e_n
\|_\infty<\epsilon $ and $\ \tau(\1-e_n )<\delta \text { for } n\geq
N_0$.
\end{defn}

\begin{rem}
We will use terms \textbf{converges in measure} \index{converges in measure} and
\textbf{converges stochastically} \index{converges stochastically} interchangeably.
\end{rem}
\begin{defn}
Let $x$ be a measurable operator from $S(M)$ and $t > 0$. The
\textbf{t-th singular number \index{singular number} of $x$} is defined as
\begin{equation}
\label{e:MCONV}
 \mu_t(x) =
inf\{||xe|| \text{ where } e \text{ is a projection in } P(M) \text{
with } \tau(\1 - e) < t \}.
\end{equation}
\end{defn}

\begin{rem}
Note that \textbf{measure topology} \index{measure topology} is defined in Fack and Kosaki's
\cite{FackKos86} as linear topology with fundamental system of
neighborhoods around $0$ given by $V(\epsilon, \delta) = \{x\in S(M)
$ such that there exists a projection $e(x,\epsilon,\delta)$ with
$\|xe\|<\epsilon$ and $\tau(\1-e)<\delta \}$.
\end{rem}

\begin{defn}
\par Denote by \textbf{
$\lambda_t(x)$ \index{$\lambda_t(x)$} the distribution function of $x$} defined as
\begin{equation}
\lambda_t(x)=\tau(E_{(t,\infty)}(|x|)), \ \ t\geq 0,
\end{equation}
here $E_{(t,\infty)}(|x|))$ is a spectral projection of $x$
corresponding to interval $(t,\infty)$.
\end{defn}

\begin{rem}
\label{r:RIGHTCONT} For the measurable operator $x$, we have
$\lambda_t(x)< \infty$ for large enough $t$ and
$\lim_{t\to\infty}\lambda_t(x)=0$. Moreover, the map $\R  \ni t  \to
\lambda_t(x)$, is non-increasing and continuous from the right
(because $\tau$ is normal and $\{|x|>t_n\}\uparrow \{|x|>t\}$ (and
hence in strong operator topology) as $t_n\downarrow t$). The distribution
$\lambda_t(x)$ is a non-commutative analogue of the distribution
function in classical analysis, (see. \cite{FackKos86} p. 272 or
\cite{SteinWeiss71}).
\end{rem}
We would need the following statement about properties of the $\mu_t(x)$
(see for example proposition 2.4 \cite{Yeadon76}, or lemma 2.5
\cite{FackKos86}):
\begin{lem}
Let $x,y\in S(M)$ be measurable operators.
\begin{itemize}
  \item[i)] Map $\ \R \ni t \to \mu_t(x)$ is non-decreasing and
  continuous from the right.

  Moreover $lim_{t\downarrow  0}\mu_t(x)=\|x\|_\infty\in [0,\infty]$,
  \item[ii)] $\mu_t(x)=\mu_t(|x|)=\mu_t(x^*)$ and $\mu_t(\alpha x) =
  |\alpha| \mu_t(x)$ for $\alpha\in \mathbb{C},\ t>0$,
  \item[iii)] $\mu_t(x)\leq \mu_t(y)$ for $\ 0\leq x\leq y,
  \ t>0$,
  \item[iv)] $\mu_{t+s}(x+y)\leq \mu_t(x)+\mu_s(y)$ for $t,s >0$,
  \item[v)] $\mu_t(yxz)\leq \|y\|_\infty\|z\|_\infty\ \mu_t(x)$, for
  $y,z\in M,\ \ t>0$,
  \item[vi)] $\mu_{t+s}(yx)\leq \mu_t(x)\mu_s(y)$ for $\
  t,s>0$.
\end{itemize}
\end{lem}

\section{\protect\smallskip Stochastic Banach Principle \index{Stochastic Banach Principle}}
We start the section with the description of some conditions
equivalent to stochastic convergence (cmp. with lemma 3.1
\cite{FackKos86}).
\begin{lem}
Let $M, \ \tau$ be as before. Consider following conditions:
\begin{itemize}
  \item[i)]  Sequence $\{x_n\}_{n=1}^\infty$ converges to $0$ in
  measure,
  \item[ii)]For every $\ \epsilon>0, \delta>0$ there
    exist a positive real $\ 0<\delta'<\delta$ and an integer $N_0$ such that for $n\geq N_0$
    $$\mu_{\delta'}(x_n)<\epsilon,$$
    \item[iii)]  For every $\ \epsilon>0, \delta>0, p\in
P(M)$ with $\tau(p)<\infty$ there exists an integer $N_0$ and a
sequence of projections $
      \{e'_n\}_{n\geq N_0} \subset P(M),\ e'_n\leq p $ such that
      $$
     \| x_n\, e'_n \|_\infty<\epsilon \text{ and } \ \tau(p-e'_n )<\delta \text { for } n\geq
     N_0.
    $$
\end{itemize}
The following relations take place: $ i)\Leftrightarrow ii) \Rightarrow
iii)$. If $\tau$ is finite then $iii)\Rightarrow i)$.
\begin{proof}
Implication $ ii)\Rightarrow i)$  follows from the fact that
condition $\mu_{\delta'}(x_n)<\epsilon$ implies for the sequence of
projections $\{e_n\}_{n=1}^\infty,$ holds
$\|x_n\,e_n\|\leq 2\epsilon $ and $\tau(\1-e_n)\leq\delta'$.
Implication $ i)\Rightarrow ii)$ follows from the definition of
measure convergence \ref{d:STOCH}.

Implication $i), ii)\Rightarrow iii)$ follows from the inequality

$\tau(p - p\wedge q) = \tau(p\vee q - q)\leq \tau(\1-q)$, hence
sequence $\{e'_n=e_n\wedge p\}_{n=1}^\infty$ satisfies
iii) (here projections $e_n$ are defined in the proof $ii)\Rightarrow i))$.

The case when $\tau$ is finite follows immediately since
$\tau(\1)<\infty$.
\end{proof}
\end{lem}

We need the following technical statement which is interesting by
itself:
\begin{lem}
Let $x,y$ be self-adjoint measurable operators from $S(M)$, $t,s$ be
 positive real. Then
\begin{equation}
\label{e:SPECTR} \lambda_{t+s}(x+y)\leq\lambda_t(x)+\lambda_s(y)
\end{equation}
\begin{proof}
Indeed,
\begin{equation}
\begin{split}
\label{e:LONGNORM} &\|
|(x+y)\,|(\1-\{|x|>t\})\wedge(\1-\{|y|>s\})\|=\\&
\|(x+y)\,(\1-\{|x|>t\})\wedge(\1-\{|y|>s\})\|\leq\\
&\leq\|x\,(\1-\{|x|>t\})\wedge(\1-\{|y|>s\})\|+\|y\,(\1-\{|x|>t\})\wedge(\1-\{|y|>s\})\|=\\
&=\| |x|\,(\1-\{|x|>t\})\wedge(\1-\{|y|>s\})\|+\|
|y|\,(\1-\{|x|>t\})\wedge(\1-\{|y|>s\})\|\leq\\
&\leq \| |x|\,(\1-\{|x|>t\})\|+\| |y|\,(\1-\{|y|>s\})\|\leq\\&t+s
.\end{split}
\end{equation}
Here the first and the second equality follows from the equality
$\||z|\,u_z^*u_z\,|z|\|=\||z|^2\|=\|z^*\,z\|$, where $z\in M_h$,
$u_z$ is a partial isometry from $M$ such that $z=u_z\,|z|$, and
\begin{equation}
\label{e:SPECTRAL} u^*_zu_z=l(z),\ u_zu^*_z=r(z) ,
\end{equation}
where $l(z)(r(z))$ is a left (right) support of $z$. Inequality
\ref{e:LONGNORM} means that
\begin{equation}
\mu_{\lambda_t(x)+\lambda_s(y)}(x+y)\leq t+s.
\end{equation}
Let $\xi$ be a vector from Hilbert space $\Hs$ and suppose that
\begin{equation}
\xi\in\ \{|x+y|>s+t\}\Hs\cap (\1-\{|x|>t\})\wedge(\1-\{|y|>s\})\Hs.
\end{equation}
Then
\begin{equation}
\label{e:HSINEQ} ((t+s)\|\xi\|)^2 < (|x+y|\xi,|x+y|\xi) = ((x+y)\xi,(x+y)\xi)
\leq ((t+s)\|\xi\|)^2,
\end{equation}
here the first inequality follows from inclusion $\xi\in\ \{|x+y|>s+t\}\Hs$,
the equality follows from the
spectral decomposition \ref{e:SPECTRAL}, the second inequality follows
from inclusion $\xi\in\
(\1-\{|x|>t\})\wedge(\1-\{|y|>s\})\Hs$.

Inequality \ref{e:HSINEQ} implies that $\|\xi\|=0$ or, in other words,
\begin{equation*}
\{|x+y|>s+t \} \wedge ((\1-\{|x|>t\})\wedge(\1-\{|y|>s\}))=0.
\end{equation*}
Hence,
\begin{equation}
\begin{split}
&\{|x+y|>t+s\}=\{|x+y|>t+s\}-\\&
\{|x+y|>t+s\}\wedge((\1-\{|x|>t\})\wedge(\1-\{|y|>s\}))\thicksim\\
&\thicksim\{|x+y|>t+s\}\vee((\1-\{|x|>t\})\wedge(\1-\{|y|>s\}))-\\
&((\1-\{|x|>t\})\wedge(\1-\{|y|>s\}))\leq\\
&\leq\1-((\1-\{|x|>t\})\wedge(\1-\{|y|>s\}))=\{|x|>t\}\vee\{|y|>s\}.
\end{split}
\end{equation}
Here $\thicksim$ means projection equivalence. Since trace $\tau$ is
invariant on equivalent projections,
\begin{equation}
\tau(\{|x+y|>t+s\})\leq \tau(\{|x|>t\}\vee\{|y|>s\})\leq
\tau(\{|x|>t\})+ \tau(\{|y|>s\})
\end{equation}

and, hence, the inequality (\ref{e:SPECTR}) takes place.
\end{proof}
\end{lem}

\begin{thm}
\label{t:STOCBANAH} Let $( B, \|.\|)$ be a Banach space. Let $\Sigma
= \{A_n, n \in \N\}$ be a set of linear operators $A_n: B \to
S(M)$.
\begin{itemize}
  \item [i)] Suppose that there exists a function $C(\lambda):\R_+\to\R_+$ with
  $lim_{\lambda\to\infty}C(\lambda)=0$, and such that
\begin{equation}
\label{e:UNIFBOUND}
\sup_{n\in\N}\tau(\{|A_n(b)|>\lambda\,\|b\|\})\leq C(\lambda)
\end{equation}
holds for every $b\in B,\ \lambda\in\R_+$.

Then the subset $\tilde{B}$ of $B$ where $A_n(b)$ converges in measure
(stochastically) is closed in $B$.
   \item [ii)] Conversely, if $A_n$ is a set of continuous in
measure maps from $B$ into $S(M)$ and for each $b\in B, \lambda\in
R_+$
\begin{equation}
\lim_{\lambda\to\infty}\sup_{n\in\N}\tau(\{|A_n(b)|>\lambda\})=0,
\end{equation}
then there exists a function $C(\lambda):\R_+\to\R_+$ with
$\lim_{\lambda\to\infty}C(\lambda)=0$, and
\begin{equation}
\sup_{n\in\N}\tau(\{|A_n(b)|>\lambda\,\|b\|\})\leq C(\lambda)
\end{equation}
\end{itemize}

Part i) of the theorem \ref{t:STOCBANAH} means that under the
condition of the linear uniform boundness (\ref{e:UNIFBOUND}), the set of
the stochastical convergence is closed.

Part ii) of the theorem \ref{t:STOCBANAH} means that if the set of
uniform boundness  is closed, then linear uniform boundness takes
place.

Note that even the condition in the part ii) looks more restrictive,
it is similar in nature to the condition of part i), since we can
restrict everything to the closed liner subspace $B_1$
 of Banach space $B$ ($B_1$ is also Banach space).

\begin{proof}

Part i) We first show that condition \ref{e:UNIFBOUND} implies
continuity of the set $\Sigma$ of operators . Let $B\supset
\{b_k\}_{k=1}^\infty$ be a sequence in $B$ converging to $b\in
B$. Then for $\lambda, \epsilon\in\R_+$ with $2\lambda\,\sup_{k\geq
n}\|b-b_k\|<\epsilon$
\begin{equation}
\begin{split}
\tau(\{|A_n(b_k)-A_n(b)|>\epsilon\})\leq
\tau(\{|A_n(b_k-b)|>\lambda\,\|b_k-b\|\})\leq\\
C(\lambda\,\|b_k-b\|^{-1})\stackrel{b_k\rightarrow b}\longrightarrow 0,
\end{split}
\end{equation}
for $k\geq n $, and, hence $A_n$ continuous. Note that the inequality follows from the fact that
right part of \ref{e:UNIFBOUND} does not depend on the norm of $b$.

There exists a subsequence $b_{k_j}$ of $b_k$ such that sequence
$x_j=lim_{n\rightarrow\infty}A_n(b_{k_j})$ converges stochastically.
To show this we choose a sequence of $\{k_i'\}_{i=1}^\infty$ base on the inequality
\ref{e:UNIFBOUND} in such a manner that
\begin{equation}
\label{e:ESTBK} \tau(\{|A_n(b_{k_j'}-b_{k_{j+l}'})|>2^{-j}\})\leq
2^{-j} \mbox{ for all n }
\end{equation}
and
\begin{equation}
\label{e:ESTBKL} \tau(\{|A_n(b_{k_j'}-b)|>2^{-j}\})\leq 2^{-j}
\mbox{ for all n }
\end{equation}

This may be done since $b_n\stackrel{n\rightarrow
\infty}\longrightarrow b$, and
$C(\lambda)\stackrel{\lambda\rightarrow \infty}\longrightarrow 0$.
It is sufficient to choose a sequence of $\{\lambda_j\}_{j=1}^\infty$ in such a way
that $C(\lambda_j)<2^{-j}$ and
$\|b_{k_j}-b\|<\lambda_j^{-1}\,2^{-2j}$.

Choose $n_j$ in such a manner that for $N>n_j$ holds
\begin{equation}
\label{e:XLIMITK} \tau(\{|A_N(b_{k_j})-x_j|>2^{-j}\})<2^{-j}.
\end{equation}

This is possible since $A_n(b_{k_j})$ converges stochastically to
$x_j$.

Then for $j,i \in\N$ and $n>n_{i+j}$
\begin{equation}
\begin{split}
\tau(\{|x_j-x_{j+i}|>3\,\cdot 2^{-j}\})= &\tau(\{|(x_j-A_n(b_{k_j}))+\\
(A_n(b_{k_j})-A_n(b_{k_{j+i}}))+&(A_n(b_{k_{j+i}})-x_{j+i})|>3\,\cdot 2^{-j}\})\leq\\
\tau(\{|A_n(b_{k_j})-x_j|>2^{-j}\})+&\tau(\{|A_n(b_{k_j})-A_n(b_{k_{j+i}})|>2^{-j}\})+\\
\tau(\{|A_n(b_{k_{j+i}})-x_{j+i}|>&2^{-(j+i)}\})\leq3\cdot 2^{-j}.
\end{split}
\end{equation}
Here the first inequality follows from the \ref{e:SPECTR}.

Denote the stochastic limit of $\{x_j\}_{j=1}^\infty$  by $x_0$. If
necessary by taking a subsequence of $\{x_j\}$ and reindexing, we suppose
that
\begin{equation}
\label{e:XLIMIT} \tau(\{|x_j-x_0|>2^{-j}\})\leq 2^{-j}.
\end{equation}
Sequence $\{A_n(b)\}_{n=1}^\infty$ converges to $x_0$
stochastically. Indeed, for $n>n_j$ holds the following inequality holds
\begin{equation}
\begin{split}
\tau(\{|A_n(b)-x_0|>3\,\cdot 2^{-j}\})= &\tau(\{|(A_n(b) -
A_n(b_{k_j}))+ \\
(A_n(b_{k_j})-x_j)+ (x_j-x_0)|>3\,\cdot 2^{-j}\}) \leq &
\tau(\{|A_n(b) - A_n(b_{k_j})|>2^{-j}\})+\\
\tau(\{|A_n(b_{k_j})-x_j|>2^{-j}\})+ &
\tau(\{|x_j-x_0|>2^{-j}\})\leq 3 \cdot 2^{-j}.
\end{split}
\end{equation}
Here the first inequality follows from \ref{e:SPECTR}, and the
second inequality follows by noting that the first part follows from
\ref{e:ESTBKL}, the second part follows from \ref{e:XLIMITK} and
choice of $n$, and the third part follows from \ref{e:XLIMIT}.

Part i) is established.

Part ii) Suppose that for every $b\in B$ and $\lambda\in \R_+$ holds
\begin{equation}
\label{e:UNIFCONTO}
\sup_n\tau(\{|A_n(b)|>\lambda\})\stackrel{\lambda\rightarrow\infty}\longrightarrow
0.
\end{equation}
For fixed $\epsilon>0$ and $\lambda\in \N$ define $B_\lambda =
\{b\in B | \sup_n \tau(\{|A_n(b)|>\lambda\})\leq\epsilon\}$. Then
from \ref{e:UNIFCONTO} it follows that
\begin{equation}
B=\bigcup_{\lambda\in \N} B_\lambda
\end{equation}
Let $B_{\lambda,k}$ be a set defined as $\{b\in B |
\sup_{n\geq k}\tau(\{|A_n(b)|>\lambda\})\leq\epsilon\}$. Then
\begin{equation}
\label{e:BN} B_\lambda=\bigcap_{k\in \N} B_{\lambda,k}
\end{equation}
Sets ${B_{\lambda,k}}$ are closed. Indeed, let
$B_{\lambda,k}\supset\{b_j\}_{j=1}^\infty$ converges to $b\in B$.
Then
\begin{equation}
\label{e:N+H}
\begin{split}
\tau(\{|A_n(b)|>\lambda+\gamma\})=&\tau(\{|A_n(b_k)-(A_n(b_j)-A_n(b))|>\lambda+\gamma\})\leq\\
\tau(\{|A_n(b_j)|>\lambda\}+&\tau(\{|(A_n(b_j)-A_n(b))|>\gamma\})\leq\epsilon
\end{split}
\end{equation}
Here the first inequality follows from \ref{e:SPECTR}, the
first estimate follows from the definition of ${B_{\lambda,k}}$ and the
second estimate become valid for sufficiently large $j$,
and follows from the free choice of $b_j$ and continuity of
$A_n$ in measure.

Since $\lambda_t(x)$ is continuous from the right
(\ref{r:RIGHTCONT}), then
\begin{equation}
\tau(\{|A_n(b)|>\lambda\})=\lim_{m\rightarrow
\infty}\tau(\{|A_n(b)|>\lambda+\gamma_m\})\leq\epsilon,
\end{equation}
where $\gamma_m\stackrel{m\rightarrow\infty}\longrightarrow 0$.
Hence, $b\in B_{\lambda,k}$, or $ B_{\lambda,k}$ is closed. Set
$B_\lambda$ is closed as an intersection of closed sets (\ref{e:BN}).

It follows from the Baire category principle \index{Baire category principle} that there exists
$\lambda$ such that set $B_\lambda$ has non empty interior. Let
$B(b_0,r)=\{b\in B | \|b-b_0\|\leq r\}$ be contained in the
$B_\lambda$.

Then
\begin{equation}
\tau(\{|A_n(b)|>\lambda\})\leq\epsilon \mbox{ for every } b\in
B(b_0,r).
\end{equation}
Moreover, for $b=b_0-r\cdot c\in B(b_0,r)$ with $c\in B,\
\|c\|\leq 1$, holds
\begin{equation}
\label{e:RC}
\begin{split}
\tau(\{|A_n(r\cdot c)|>2\cdot\lambda\})= \tau(\{|A_n(r\cdot c -&
b_0) +A_n(b_0)|>2\cdot\lambda\})\leq \\
\tau(\{|A_n(r\cdot c - b_0)|>\lambda\})+ &
\tau(\{|A_n(b_0)|>\lambda\}) \leq 2\cdot \epsilon.
\end{split}
\end{equation}
Let $\gamma\geq 2\cdot\lambda/r$. From \ref{e:RC} it follows that
$\tau(\{|A_n(c)|>\gamma\})\leq 2\cdot\epsilon$, for every $c\in B,\
\|c\|\leq 1$.

Let $C(\gamma)=\sup_{c\in B,\ \|c\|\leq 1}
\tau(\{|A_n(c)|>\gamma\})\leq 2\cdot\epsilon$. Free choice of
$\epsilon$ implies that
\begin{equation}
\lim_{\gamma\rightarrow\infty}C(\gamma)=0,
\end{equation}
hence \ref{e:UNIFBOUND} is valid.
\end{proof}
\end{thm}
For the application of theorem \ref{t:STOCBANAH} it is
convenient to combine both parts i) and ii).
\begin{thm}
\label{t:STOCBANAHCOMBI} Let $( B, \|.\|)$ be a Banach space. Let
$A_n$ be a set of continuous in measure linear maps from $B$ into
$S(M)$, let $\lambda\in R_+$, and for each $b\in B $ holds
\begin{equation}
\label{e:UNIFORMBOUND}
\lim_{\lambda\to\infty}\sup_{n\in\N}\tau(\{|A_n(b)|>\lambda\})=0.
\end{equation}
Then subset $\tilde{B}$ of $B$ where $A_n(b)$ converges in measure
(stochastically) is closed in $B$.
\begin{proof}
Follows immediately from applying consecutively  Theorem
\ref{t:STOCBANAH} part ii) then part i).
\end{proof}
\end{thm}

Let $e$ be a projection in $M$, let $M_e$ be von Neumann algebra
consisting of operators of form $exe,\ x\in M$. If $\tau$ is a
semifinite normal faithful trace on $M$ then $\tau_e=\tau|_{M_e}$ is a
semifinite (possibly finite) faithful normal trace on $M_e$. Indeed,
tracial property, semifiniteness, normalness and faithfulness of
$\tau_e$ follows directly from similar properties of $\tau$. Space
$S(M_e,\tau_e)$ is isomorphic to the $S(M,\tau)_e$ since both these
spaces are closures of the $(M_{\tau-finite
support})_e=(M_e)_{\tau_e-finite support}$.

\begin{prop}
\label{p:1} Let $B_n$ be a sequence of continuous in measure
operators on $S(M,\tau)$. Let $e_i\in P(M),\ i=1,2,\ \1=e_1+e_2$ be
projections in $M$. Suppose that relation $e_i(B_n(x))=B_n(x_{e_i})=(B_n(x))e_i$ holds
for every $n\in\N$ and $x\in S(M,\tau)$, or, in other words, $e_i$
commutes with $B_n$. Suppose also following relations hold
\begin{equation}
\label{e:UNIFORMBOUND1}
\lim_{\lambda\to\infty}\sup_{n\in\N}\tau(\{|B_n(x_{e_i})|>\lambda\})=0,
\end{equation}
for $i=1,2$ and every $x\in S(M,\tau)$. Then the following equality
is valid:
\begin{equation}
\label{e:UNIFORMBOUND2}
\lim_{\lambda\to\infty}\sup_{n\in\N}\tau(\{|B_n(x)|>\lambda\})=0.
\end{equation}
\begin{proof}
The following relations are valid:
\begin{equation}
\label{e:UNIFORMBOUNDa}
\tau(\{|B_n(x_{e_i})|>\lambda\})=\tau(e_i\{|B_n(x)|>\lambda\}).
\end{equation}
Indeed, since for $x\in S_h(M)$ ($S_h(M)$ is a set of all
self-adjoint operators in $S(M)$) limit
$\tau(\{|x|>\lambda\})\stackrel{\lambda\rightarrow\infty}\longrightarrow
0$, hence for a sequence of polynomial $P_j(y)$ in $\R$ converging
to $\chi_{\{|y|>\lambda\}}(y)$ pointwise, sequence $P_j(x)$
converges to $\chi_{\{|x|>\lambda\}}(x)$
 stochastically. Then by \cite{FackKos86} Proposition 3.2
\begin{equation}
\begin{split}
\tau(\{|B_n(x_{e_i})|>\lambda\}) = \lim_j &\tau(P_j(B_n(x_{e_i})))=\lim_j \tau(P_j(B_n(e_ixe_i)))=\\
\lim_j \tau(P_j(e_iB_n(x)e_i))=&\lim_j
\tau(e_iP_j(B_n(x)))=\tau(e_i\{|B_n(x)|>\lambda\}).
\end{split}
\end{equation}
Statement \ref{e:UNIFORMBOUND2} follows now from the fact that (it
follows from $B_n$ commutes with $e_i$ and \ref{e:SPECTR})
\begin{equation}
\begin{split}
\tau(\{|B_n(x)|>\lambda_1+\lambda_2\})= &\tau(\{|(e_1+e_2)B_n(x)(e_1+e_2)|>\lambda_1+\lambda_2\})=\\
\tau(\{|(e_1B_n(x)e_1 &+e_2B_n(x)e_2|>\lambda_1+\lambda_2\})\leq\\
 \tau(\{|(B_n(x_{e_1})|& >\lambda_1\})+ \tau(\{|B_n(x_{e_2})|>\lambda_2\}).
\end{split}
\end{equation}
\end{proof}
\end{prop}

\begin{rem}
\label{r:1}We are going to use \ref{e:UNIFORMBOUND} in the next section
when dealing with stochastic ergodic theorem, since under the conditions of the stochastic ergodic theorem
estimate \ref{e:UNIFORMBOUND} has place.
\end{rem}

\section{\protect\smallskip Stochastic ergodic theorems}

In this section we establish stochastic convergence of the bounded
Besicovitch sequences \index{bounded
Besicovitch sequence}, and show stochastic ergodic theorems for
uniform subsequences.

In this section we use following assumptions: $M$ is a von Neumann
algebra with faithful normal tracial state $\tau$, and $\alpha$ is
an $*$-automorphism of algebra $M$. Denote by $A_n
(x)=\frac1{n}\sum_{l=1}^{n-1}\alpha^l(x),\ \mbox {for } x\in M$.
Define $\alpha'$ as a linear map on $L_1(M,\tau)$ satisfying
$\tau(x\cdot\alpha(y))=\tau(\alpha'(x)y)$ for $x\in L_1(M,\tau),
y\in M$, and $A'_n (x)=\frac1{n}\sum_{l=1}^{n-1}\alpha^{'l}(x),\
\mbox {for } x\in L_1(M,\tau)$.

Let us recall some definitions from Grabarnik and Katz
\cite{GrabKatz1995} and Chilin Litvinov and Skalski \cite{Chil}.
\begin{defn}
A positive operator $h\in M_+$ is called \textbf{ weakly wandering}  \index{weakly wandering}
if
\begin{equation}
\|A_n(h)\|_\infty\stackrel{n\rightarrow\infty} \longrightarrow 0
\end{equation}
\end{defn}
The following definition is due to Ryll-Nardzewski \cite{Ryll}.
\begin{defn}
Let $\C_1$ denote the unit circle in $\mathbb{C}$. A
\textbf{trigonometric polynomial} \index{trigonometric polynomial} is a map $P_k(n):\, \N\mapsto \C$,
where $P_k(n)=\sum_{j=0}^{k-1} b_j\cdot\lambda_j^n$ for
$\{\lambda_j\}_{j=0}^{k-1}\subset\C_1$.
\end{defn}

Bounded Besicovitch sequences are bounded sequences from the
$l_1$-average closure of the trigonometric polynomials.

More precisely,
\begin{defn}
A sequence ${\beta_n}$ of complex numbers is called a
\textbf{Bounded Besicovitch sequence} (BB-sequence) if
\begin{itemize}
  \item [(i)] $|\beta_n|\leq C < \infty $ for every $n\in\N$ and
  \item [(ii)] For every $\epsilon > 0$, there exists a trigonometric polynomial $P_k$ such that
  \begin{equation}
    \label{e:BESICEST}
     \limsup _n \frac {1}{n}\overset{n-1} {\underset{ j=1 } \sum}\; |\beta
        _{j}- P_k(j)|< \epsilon
\end{equation}
\end{itemize}
\end{defn}

Let $\mu$ be the normalized Lebesgue measure (Radon measure) on
$\C_1$. Let $\tilde{M}$ be the von Neumann algebra of all
essentially bounded ultra-weakly measurable functions $f: (\C_1,
\mu) \rightarrow M$. Algebra $\tilde{M}$ is isomorphic to
$L_\infty(\C_1,\mu)\overline {\bigotimes} M$ -which is a $W^*$
tensor product of $L_\infty(\C_1,\mu)$ and $M$, $\tilde{M}$ is a
dual to the space $L_1(\C_1,\mu)\overline {\bigotimes} M_*$ (
for definition of $W^*$ tensor product and form of the
predual space of the $W^*$ tensor product  see for example Takesaki,
\cite{Tak79}, Theorem IV.7.17).
The space $L_1(\C_1,\mu)\overline {\bigotimes} M_*$ maybe considered as a set
of $L_1$ functions on $(\C_1,\mu)$ with values in $M_*$. Algebra
$\tilde{M}$  has a natural trace $\tilde{\tau} (f) = \int_{\C_1}
\tau (f(z))d\mu (z)$, and $\tilde{M}_*$ is isomorphic to
$L_1(\tilde{M},\tilde{\tau})$.

Let $\sigma$ be an automorphism of $(\C_1, \mu)$ as a Lebesgue space
with measure. We define automorphism $\alpha\bigotimes\sigma$ of
$(\tilde{M},\tilde{\tau})$ as a closure of the linear extension of
automorphism acting on $(\tilde{M},\tilde{\tau})\ni x(z)$ as
$\alpha\bigotimes\sigma(x(z))=\alpha(x(\sigma(z)))$.
\begin{ex}
\label{ex:1} An example of such an automorphism is
$\tilde\alpha_\lambda(x(z))=\alpha(x(\lambda\cdot z)))$, for
$\lambda\in\C_1$.

In this case
\begin{equation}
A_n(x)=\frac1{n}\sum_{l=1}^{n-1}\tilde{\alpha_\lambda}^l(x) =
\frac1{n}\sum_{l=1}^{n-1}\alpha^l(x(\lambda^l\cdot z)).
\end{equation}

In particularly, if $x(z)\equiv z\cdot x $ for $x\in M$ then
\begin{equation}
A_n(x\cdot z
)=z\cdot\frac1{n}\sum_{l=1}^{n-1}\lambda^l\cdot\alpha^l(x).
\end{equation}
\end{ex}
The following lemma connects stochastic convergence in
$L_1(\tilde{M},\tilde{\tau})$ with pointwise convergence on $C_1$
and stochastic convergence in $M$ (cmp. with \cite{Chil}).

\begin{lem}
\label{t:POINWISESTOCH}
\begin{itemize}
  \item[i)] If $ L_1(\tilde{M},\tilde{\tau}) \ni x_n \stackrel{n\rightarrow\infty}\longrightarrow x_0
\in L_1(\tilde{M},\tilde{\tau})$ b.a.u. , then $
x_n(z)\stackrel{n\rightarrow\infty}\longrightarrow x_0(z)$
stochastically for almost every $ z \in \C_1$
\item[ii)] Suppose that $h$ is a weakly wandering operator with support
$supp(h)=\1$ for sequence $A_n$. Then $A_n'(x)$ converges to $0$
stochastically.
  \item[iii)] Let algebra $\mathcal{N}=(M,\tau)\overline{\otimes}{L_\infty(X,\mu)}$, (here $X$
 is a separable Hausdorff compact set, and $\mu $ is Lebesgue measure),
$\alpha$ is an automorphism of $M$,
and $\sigma$ is an automorphism of ${L_\infty(X,\mu)}$. Then
$\alpha\bigotimes\sigma$ is an automorphism of $\mathcal{N}$.
Suppose that $h$ is a weakly wandering operator with support
$supp(h)=\1$ for sequence $A_n$ corresponding to automorphism
$\alpha\bigotimes\sigma$. Then $A_n'(x(z))$ converges to $0$
stochastically for almost every $z\in\C_1$.
\end{itemize}
\end{lem}

\begin{proof}
Part i) follows from \cite{Chil}, Lemma 4.1 which states that under
the hypothesis of part i) b.a.u. convergence of $x_n(z) $ to $x_0(z)$ for
almost every z in $\C_1$, (hence double side stochastic
convergence), and the fact that double side stochastic convergence
is equivalent to (one sided) stochastic convergence (see \cite{Chil},
Theorem 2.2).

Part ii) We suppose that $x \in L_1(M,\tau)_+$ and $A_n'(x)$ is a
sequence satisfying
\begin{equation}
\tau(A_n'(x)h)\rightarrow 0 \mbox { for } n\rightarrow\infty.
\end{equation}

The following inequality is valid:
\begin{equation}
\label{e:MULTIPLSPECTR} t s \cdot\tau(\{A_n'(x)>
t\}\wedge\{h>s\})\leq \tau(A_n'(x)h).
\end{equation}
Indeed, for projections $e_1,e_2\in P(M)$, we have $e_1e_2e_1\geq
e_1\wedge e_2$. To see that note that since $e_1\wedge e_2$ commutes with $e_1, e_2$,
we have $(\1-e_1\wedge e_2)e_1e_2e_1(e_1\wedge e_2)=0$, and, hence
$e_1e_2e_1=(\1-e_1\wedge e_2)e_1e_2e_1(\1-e_1\wedge e_2)+$ $
(e_1\wedge e_2)e_1e_2e_1(e_1\wedge e_2)$ $=(\1-e_1\wedge
e_2)e_1e_2e_1(\1-e_1\wedge e_2)+(e_1\wedge e_2)$.

Then,
\begin{equation}
\begin{split}
ts\cdot \tau(\{A_n'(x)> t\}\wedge\{h>s\}) \leq& t \tau(\{A_n'(x)> t\} s \{h>s\} \{A_n'(x)> t\})\leq\\
t \tau(\{A_n'(x)> t\} h \{A_n'(x)> t\})=& t \tau( \{A_n'(x)> t\} h)\leq\\
& \tau(A_n'(x) h).
\end{split}
\end{equation}
Hence, \ref{e:MULTIPLSPECTR} is valid.

Furthermore,
\begin{equation}
\label{e:MAINIIINEQ} \tau(\{A_n'(x)> t\})\leq\frac 1{ts}
\tau(A_n'(x) h) + \tau(\1-\{h>s\}).
\end{equation}
The latter inequality follows from \ref{e:MULTIPLSPECTR}, and the
fact that $\tau(e_1)\leq\tau(e_1\wedge e_2)+ \tau(\1-e_2)$. Indeed,
\begin{equation}
\begin{split}
\tau(e_1-e_1\wedge e_2)=\tau((\1-e_1\wedge e_2)e_1(\1-e_1\wedge
e_2)) =& \tau(e_1(\1-e_1\wedge e_2)e_1)
\leq\\
 \tau(e_1(\1-e_2)e_1)=& \tau(e_1(\1-e_2))\leq \tau(\1-e_2).
\end{split}
\end{equation}
Hence \ref{e:MAINIIINEQ} is valid.

Note that inequality \ref{e:MAINIIINEQ} with the fact that
$\tau(A_n'(x) h)\stackrel{n\rightarrow\infty} \longrightarrow 0$
implies that $\sup_{n\in\N}\tau(\{|A_n(b)|>\lambda\,\|b\|\})\leq
C(\lambda)$. Indeed, sequence $\{\tau(A_n'(x) h)\}_{n=1}^\infty$ is
bounded by constant $C_0$ as a converging sequence. Choose a
monotonically decreasing sequence of $\{s_j\}_{j=1}^\infty\subset \R_+$ such
that $\tau(\1-\{h>s_j\})<2^{-j}$, and $t_j=2^{j}s_j^{-1}$. Then
\begin{equation}
\tau(\{A_n'(x)> t_j\})\leq \frac 1{t_js_j} C_0+2^{-j}=(C_0+1)2^{-j}.
\end{equation}
Hence the condition of theorem \ref{t:STOCBANAHCOMBI} is satisfied.
From the theorem follows stochastic convergence of the $A'_n(x)$
since for the dense subset in $L_1(M,\tau)$ of view $x-A_k'(x)+\widetilde{x}$
(here $\widetilde{x}\in M$ is an $\alpha'$-invariant element, see \cite{Kreng},
Theorem 1.5 (iii) on page 273 ) for
$x\in M\cap L_1(M, \tau)$ convergence is in $L_1$, and hence
stochastically.

Part iii). The proof follows line of the proof for ii). We provide only
necessary modifications. Let $E_1$ be a conditional expectation
with respect to trace $\tau\otimes\mu$ of $(
M,\tau)\overline{\otimes}{L_\infty(X,\mu)}$ onto $(
M,\tau)\overline{\otimes}{Const(X,\mu)}$, and $E_2$ be a conditional
expectation with respect to trace $\tau\otimes\mu$
of $( M,\tau)\overline{\otimes}{L_\infty(X,\mu)}$ onto
$\C\cdot\1\overline{\otimes}{L_\infty(X,\mu)}$, (for definition of
conditional expectation trace with respect to $\tau\otimes\mu$
and its existence see \cite{Tak79}). Due to
the form of the $\alpha\bigotimes\sigma$, both $E_j$'s commute with
$A_n$, for $j=1,2$.

Since
\begin{equation}
\label{e:WWH} \|A_n(h)\|_\infty\geq \|E_1 A_n(h)\|_\infty= \|
A_n(E_1 h)\|_\infty,
\end{equation}
 and $supp(h)\leq supp(E_1 h)$ is valid, it follows that $supp(E_1
h)=\1$.  Indeed, $x\geq 0,\ x\neq 0$ implies $\tau(E_1 x)=\tau(x)>0$
hence $0<\tau((E_1 a) h)=\tau(a (E_1 h))$ and $supp(E_1 h)=\1$ for every $a\in M$.

Hence $E_1(h)$ is a weakly wandering operator.

For positive $x(z)\in L_1(M,\tau)\overline{\otimes}L_1(X,\mu)$ holds
\begin{equation}
\|x\|_1=\int_{X}\|x(z)\|_1\cdot d\mu(z),
\end{equation}
hence $\|x(z)\|_1 $ is an $L_1(X,\mu)$ function. Applying classical Hopf
inequality (see for example \cite{Kreng}, Theorem 2.1, p. 8) we get
\begin{equation}
\mu(\sup_n\{\|A_n'(x)(z)\|_1>\lambda\}) \leq \frac {Const}\lambda
\int_{X}\|x(z)\|_1\cdot d\mu(z),
\end{equation}
or, outside of a set $X_0\subset X$ of small measure the value of
$\|A_n'(x)(z)\|_1$ is uniformly bounded. Proceeding like in the part
ii) applied for every $z\in X_0$, we get stochastic converges for
every $z\in X_0$.
\end{proof}

\begin{thm}[Neveu Decomposition for the special case of tensor product of von Neumann algebras]
\label{t:NEVEU} Let algebra
$\mathcal{N}=(M,\tau)\overline{\otimes}{L_\infty(X,\mu)}$, (here $X$
 is a Hausdorff separable compact set, and $\mu $ is Lebesgue measure), $\alpha$ is an automorphism of $M$,
and $\sigma$ is an automorphism of ${L_\infty(X,\mu)}$. Then
$\tilde{\alpha}=\alpha\otimes\sigma$ is an automorphism of
$\mathcal{N}$. Suppose that in addition automorphism $\sigma$ is
ergodic. Then there exists an $\tilde{\alpha}$ invariant projection
in $\mathcal{N}$ of view $e_1 = e_{11}\otimes\1,\ e_2=\1-e_1$ with
 $e_1(z)=e_M$ for almost every $z\in X$ such that
\begin{itemize}
  \item[i)] There exists a normal state $\rho$ on $\mathcal(N)$ with $supp(\rho)=e_1$
and for almost each $z\in X$, $\rho(z)$ is invariant with respect to
automorphism $\alpha'$;
  \item[ii)] There exists a weakly wandering operator $h\in \mathcal{N}$ with $supp(h)=e_2$
and  for almost each $z\in X$, $h(z)$ is a weakly wandering operator
in $M$.
\end{itemize}
\begin{proof}
Corollary 1.1 of \cite{GrabKatz1995} implies existence of the
projection $\tilde{e_1}$ in $\mathcal{N}$ such that i) there exists
$\tilde{\alpha}'$ invariant normal state $\rho$ with support
$supp(\rho)=\tilde{e_1}$ and ii) there exists a weakly wandering
operator $h\in\mathcal{N}$ with support $\1-\tilde{e_1}$. Our goal
is to show that similar statements are valid for almost every $z\in
X$.

Since $\sigma$ is ergodic, then for every $x\in M\otimes
Const(X,\mu)$ (constant function on $X$ with values in $M$) holds
\begin{multline}
\rho(z)(x(z))=(\tilde{\alpha}'\rho(z))(x(z))= \rho(z)(\alpha(x(\sigma(z))))=\\
\rho(z)(\alpha(x(z))) = (\alpha'(\rho(z)))(x(z))
\end{multline}
,or $\rho(z)$ is $\alpha'$ invariant. Suppose that function $z\to \rho(z)$ is not
constant or $z\to \rho(z)$ is such that there exists real $r_0\in\R_+$ and
$x(z)\equiv x_0\in M_+$ with $\mu(\{z\in X |\rho(z)(x(z))\leq
r_0\})>0$ and $\mu(\{z\in X |\rho(z)(x(z))< r_0\})>0$. Since
$\sigma$ is ergodic, there exists $n\in\N$ such that
\begin{equation}
\mu(\sigma^{-n}(\{z\in X |\rho(z)(x(z))\leq r_0\})\cap \{z\in X
|\rho(z)(x(z))< r_0\})>0.
\end{equation}
Hence,
\begin{multline}
\rho(z)(x(z))=(\tilde{\alpha}^{'n}\rho(z))(x(z))=
(\alpha')^{n}(\rho(z))(x(\sigma^n(z)))=\\
\rho(z)(x(\sigma^n(z)))=\rho(\sigma^{-n}z)(x((z))),
\end{multline}
or $r_0\geq \rho(z)(x_0)=\rho(\sigma^{-n}z)(x_0)<r_0$.
Contradiction shows that function $z\to \rho(z)$ is constant.

This implies that $supp(\rho)= supp(\rho(z))=\tilde{e_1}(z)$ is
constant.

Part ii) follows directly arguments of proof \ref{e:WWH}.
\end{proof}
\end{thm}

\begin{thm}
Let algebra
$\mathcal{N}=(M,\tau)\overline{\otimes}{L_\infty(X,\mu)}$, (here $X$
 is a separable Hausdorf compact set, and $\mu $ is normalized Lebesgue measure),
 $\alpha$ is an automorphism of $M$, and $\sigma$ is an automorphism
of ${L_\infty(X,\mu)}$. Then  $\tilde{\alpha}=\alpha\otimes\sigma$
is an automorphism of $\mathcal{N}$. Suppose that in addition
automorphism $\sigma$ is ergodic. Then for almost every $z\in X$
the averages $A_n'(x(z))$ converges stochastically.
\begin{proof}
Proof of the theorem follows directly from \ref{t:NEVEU} and
\ref{t:POINWISESTOCH} applied the part where there exists weakly
wondering operator, and from the regular individual ergodic theorem
\cite{Yeadon76} applied to the part, where an invariant normal state
exists, and Proposition \ref{p:1}.
\end{proof}
\end{thm}

Now we are in a position to prove stochastic convergence of the
bounded Besicovitch sequences.

\begin{thm}[Stochastic Ergodic Theorem for bounded Besicovitch sequences]
Let $\{\beta_j\}_{j=1}^\infty$ be a bounded Besicovitch sequence.
Let $M$ be a von Neumann algebra with finite faithful normal tracial
state $\tau$. Let $\alpha$ be an automorphism of $M$. Then the sequence
$$\tilde{A_n}(x)=\frac 1{n} \sum_{j=0}^{n-1} \beta_j \alpha^{'j}(x) $$
converges stochastically for $x\in L_1(M,\tau)$.
\begin{proof}
Suppose first that bounded Besicovitch sequence
$\{\beta_j\}_{j=1}^\infty$ is a trigonometric polynomial $P_k(j)$.
Then the statement of the theorem is valid.

Indeed, choosing $\tilde{\alpha}$ as in example \ref{ex:1} we get
from theorem \ref{t:STOCBANAHCOMBI} and the fact that irrational
rotation on the $\C_1$ is ergodic (Equidistribution Kronecker-Weyl
Theorem, see for ex. \cite{Katok95} p. 146) that
\begin{equation}
\label{e:pol} A_n(x\cdot z
)=z\cdot\frac1{n}\sum_{l=1}^{n-1}\lambda^{'l}\cdot\alpha^l(x),
\end{equation}
hence
$$
\label{e:s} \frac 1{n}\sum_{l=1}^{n-1}\lambda^l\cdot\alpha{'^l}(x)
$$
converges stochastically for irrational $\lambda$.

For the rational $\lambda$ convergence follows from the fact that it
is a finite combination of averages of the $\alpha^{'m}$, where $m$
is denominator.

Taking linear combinations of terms as in  \ref{e:pol} implies the statement for
trigonometric polynomials.

Statement of the theorem is valid for the $x\in M\cap S(M)$. Indeed,
using approximation of the BB sequence by trigonometric polynomials
as in \ref{e:BESICEST} one gets for $A_n(k,x)= \frac
1{n}\sum_{l=1}^{n-1}P_k(l)\cdot\alpha{'^l}(x)$
\begin{equation}
\|\tilde{A_n}(x)-A_n(k,x)\|_\infty\leq \frac
1{n}(\sum_{l=0}^{n-1}|\beta_l-P_k(l)|)\cdot\|x\|_\infty
\end{equation}
and, hence, stochastic convergence.

Note also that for every $x\in L_1(M,\tau)$
\begin{equation}
\|\tilde{A_n}(x)-A_n(k, x)\|_1\leq \frac
1{n}(\sum_{l=0}^{n-1}|\beta_l-P_k(l)|)\cdot\|x\|_1.
\end{equation}

Hence by remark \ref{r:1} averages $\tilde{A_n}(x)$ are uniformly
bounded in the sense of \ref{e:UNIFBOUND}.

Result of the theorem follows from the Stochastic Banach Principle,
Theorem 3.2.4 and density of $M\cap S(M)$ in $L_1(M,\tau)$.
\end{proof}
\end{thm}

The following theorem is implied by the stochastic ergodic theorem for
bounded Besicovitch sequences. (cmp. \cite{Litv})

For the following definitions see for example \cite{Kreng}, p. 260.

Let $\sigma$ be a homeomorphism of a compact metric space $X$ with
metric $\varrho$ such that all powers of $\sigma^l$ are
equicontinous. Assume also that there exists $z\in X$ with dense
orbit $\sigma^l(z)$ in $X$. Then there exists a unique (hence
ergodic) $\sigma$ invariant measure $\nu$ on the $\sigma$ algebra of
Borel sets $\mathfrak{B}$. Each non-empty open set has a positive
$\nu$ measure.

A sequence $u_j$ is called \textbf{uniform} if there exists such
dynamical system $(X,\mathfrak{B},\nu,\sigma)$ and a set $Y\in
\mathfrak{B}$ with $\nu(\partial Y)=0$ and $\nu(Y)>0$ and point
$y\in X$ with $u_j= j^{th}$ entry time of orbit of $y$ into $Y$.

\begin{thm}
Let $M,\tau,\alpha$ be as in previous theorem, $\{u_j\}_{j\geq 0}$ be a
uniform sequence. Then averages
$$
\frac 1{n} \sum_{j=0}^{n-1} \alpha^{'u_j} x
$$
converge stochastically for $x \in L_1(M,\tau)$.
\begin{proof}
Follows from the previous theorem and the fact (see \cite{Ryll}) that
any uniform sequence is a bounded Besicovitch sequence.
\end{proof}
\end{thm}
\begin{rem}
Similar results remain valid for the case when $M$ is a semifinite JBW algebra with faithful
normal trace $\tau $
\end{rem}
\bibliographystyle{amsplain}

\end{document}